\documentclass[draft,a4paper]{amsart}

\usepackage{amsmath}
\usepackage{amssymb}
\usepackage{amsthm}
\usepackage[latin1]{inputenc}
\usepackage[british]{babel}

\newtheorem{Thm}{Theorem}

\theoremstyle{definition}

\theoremstyle{remark}

\newcommand{\cpx}[1]{\ensuremath{\mathbf{#1}}}          % complex
\newcommand{\digraph}[1]{G_{#1}} 
\newcommand{\morsegraph}[2]{G_{#1}^{#2}}

\newcommand{\card}[1]{|#1|}

\date{\today}

\author{Emil Sköldberg}

\address{Department of Mathematics\\ National University of Ireland, Galway\\Ireland}
\email{emil.skoldberg@nuigalway.ie}
\urladdr{www.math.su.se/\~{}emil}
\thanks{The author was supported by Marie Curie fellowship HPMD-CT-2001-00079}

\title[Homology of Heisenberg Lie algebras]{The homology of Heisenberg Lie algebras over fields of characteristic two}

\begin{document}

\maketitle

\begin{abstract}
The generating function of the Betti numbers of  the  Heisenberg Lie algebra 
over a field of characteristic $2$ is calculated using discrete Morse 
theory.
\end{abstract}

% \section{Introduction} \label{intro}

The Heisenberg Lie algebra of dimension $2n+1$, 
denoted by $\mathfrak{h}_{n}$, is the vector space with basis 
$B = \{z, x_1, \dots, x_n, y_1, \dots y_n\}$ where the only non-zero 
Lie products of basis elements are 
\[
[x_i,y_i] = -[y_i,x_i]= z.
\]
In this paper the Betti numbers of the homology of  $\mathfrak{h}_{n}$ 
over a field of characteristic $2$ 
is computed with the aid of algebraic discrete Morse theory from \cite{skoldberg:morse}.
The notation from \cite{skoldberg:morse} will be freely used.

\begin{Thm}
  The generating function of the Betti numbers of the Heisenberg Lie algebra
  over a field of characteristic $2$ is 
\[
\sum_{i \geq 0} \dim_k H_{i}(\mathfrak{h}_n) \, t^i =
 \frac{(1+t^3)(1+t)^{2n} + (t+t^2)(2t)^n}{1+t^2}
\]
\end{Thm}

When the ground field of $\mathfrak{h}_{n}$ has characteristic $0$, 
Santharoubane~\cite{santharoubane:heisenberg} has shown that
\[
 \dim_k H_{i}(\mathfrak{h}_n) = \binom{2n}{i} - \binom{2n}{i-2},
\]
(the need for the ground field to have characteristic $0$
is not explicitly mentioned).

% \section{Homology of finite-dimensional Lie algebras}

Let us first recall the construction of the Chevalley--Eilenberg complex $\cpx{V}$ of 
$\mathfrak{h}_n$, whose homology is the homology of  $\mathfrak{h}_n$:
the complex $\cpx{V}$ is given by
\[
0 \longrightarrow
\bigwedge^{2n+1} \mathfrak{h}_n
\longrightarrow \cdots
\longrightarrow \bigwedge^{i} \mathfrak{h}_n
\longrightarrow \cdots
\longrightarrow \bigwedge^{2}\mathfrak{h}_n
\longrightarrow \mathfrak{h}_n
\longrightarrow 0
\]
with the differential
\begin{equation*}
\bar{d}(w_{1} \wedge \cdots \wedge w_{n}) =
\sum_{i < j} (-1)^{i+j} [w_i,w_j] \wedge
  w_1 \wedge \cdots \wedge \widehat{w_i} \wedge  \cdots \wedge \widehat{w_j} \wedge
\cdots \wedge w_n
\end{equation*}
for $w_i \in B$.

The $p$-th homology (with trivial coefficients) of
$\mathfrak{h}_n$, can now be obtained 
as the $p$-th homology group of the complex 
$\cpx{V}$.

If $I=\{i_1,\dots,i_s\}$ is a subset of $[n]$, we will
use the notation $x_I$ for the element 
$x_{i_1}\wedge\cdots\wedge x_{i_{s}}$, (and similarly for $y_I$).

\begin{proof}
The result is proved by constructing a Morse matching $M$ on the digraph
$\digraph{\cpx{V}}$, and showing that when  $\pi$ is the projection coming from the
splitting  homotopy of $M$, we have that $\pi(\cpx{V})$ has trivial differential. 

The decomposition of the Chevalley--Eilenberg complex we will use is the obvious;
we consider the basis for $\cpx{V}$ given by
$\{z \wedge x_I \wedge y_J, x_I \wedge y_J \mid I,J \subseteq [n]\}$.

Let the matching $M$ consist of the following edges in $\digraph{\cpx{V}}$:
\[
x_I \wedge y_J \rightarrow z \wedge x_{I\setminus \{k\}} \wedge y_{J\setminus\{k\}}
\]
whenever  $\max(I^{c}\cap J^{c}) < \max(I\cap J)$ and $k = \max(I \cap J)$.

% Vilka är de omatchade elementen ?

There are now two kinds of unmatched elements:
first the elements $z \wedge x_I \wedge y_J$, with  
$\max(I^{c} \cap J^{c})  < \max(I \cap J)$, and then
the elements $x_I \wedge y_J$, with
$\max(I^{c} \cap J^{c})  > \max(I \cap J)$.

% Är matchningen acyklisk ?

When $x_I \wedge y_J \in M^{+}$,  there is exactly one element
$z \wedge x_{I'} \wedge y_{J'}$ with 
$x_I \wedge y_J \rightarrow z \wedge x_{I'} \wedge y_{J'}$ that is not
in $M^{0}$, which implies that
there can be no directed cycle in the graph $\morsegraph{\cpx{V}}{M}$. 
Together with the fact that for all edges in $\digraph{\cpx{V}}$ the corresponding 
component of the differential is an 
isomorphism, this implies that $M$ is a Morse matching.

We will now see that the differential in $\pi(\cpx{V})$ is zero.
For an element $z\wedge x_I \wedge y_J \in M^{0}$ it is obvious that
$d\pi(z\wedge x_I \wedge y_J) =  \pi d (z\wedge x_I \wedge y_J) =0$.
For  $x_I \wedge y_J \in M^{0}$ with 
$m=\max(I^{c} \cap J^{c})$ we get
that 
\[
\pi(x_I \wedge y_J) = x_I \wedge y_J + \sum_{i\in I\cap J} x_{(I\setminus\{i\})\cup\{m\}} \wedge y_{(J\setminus\{i\})\cup\{m\}} ,
\]
from which it is easily seen that $d\pi(x_I \wedge y_J) = 0$.
From \cite[Theorem 1]{skoldberg:morse} now follows that the $i$-th Betti number
is equal to the number of unmatched vertices in homological degree $i$.

For the computation of the generating function we introduce the elements  
$u_i = x_i \wedge y_i$, and we begin by counting the critical vertices 
$z \wedge x_{I} \wedge y_{J} \wedge u_{K}$ and 
$x_{I} \wedge y_{J} \wedge u_{K}$  when $I \cup J = L$ for a fixed set 
$L \subseteq [n]$.

If $L=[n]$, the critical vertices are all $z \wedge x_{I} \wedge y_{J}$ and
$x_{I} \wedge y_{J}$ and they contribute with 
$(1+t)(2t)^n$ toward the homology.

If $L \neq [n]$, then the critical vertices of the form
$z \wedge x_{I} \wedge y_{J} \wedge u_{K}$ are those with
$\max \left([n] \setminus (I \cup J)\right) \in K$ so they contribute with 
$t^3 (2t)^{\card{L}} (1 + t^2)^{n-\card{L}-1}$ toward the homology. The critical 
vertices of the form 
 $x_{I} \wedge y_{J} \wedge u_{K}$ are those with
$\max \left([n] \setminus (I \cup J)\right) \not\in K$ and thus contribute with 
$(2t)^{\card{L}} (1 + t^2)^{n-\card{L}-1}$ toward the homology.

Summing up we get 
\begin{align*}
  f(t) & =
    (1+t)(2t)^{n} + (1+t^3) \sum_{L \subset [n]}(2t)^{\card{L}}(1+t^2)^{n-\card{L}-1} \\
  & = (1+t)(2t)^{n} + (1+t^3) \sum_{i=0}^{n-1} \binom{n}{i}(2t)^i(1+t^2)^{n-i-1} \\
  & = (1+t)(2t)^n + (1+t^3)(1+t^2)^{-1} ((1+2t+t^2)^{n} - (2t)^{n}) \\
  & = \frac{(1+t)(1+t^2)(2t)^n}{1+t^2} + 
    \frac{(1+t^3)(1+t)^{2n}-(1+t^3)(2t)^n}{1+t^2} \\
    & = \frac{(1+t^3)(1+t)^{2n} + (t+t^2)(2t)^n}{1+t^2}
\end{align*}
% OBS ovanstående formel stämmer för n <= 6, kollat mot bergman
\end{proof}

\bibliographystyle{amsalpha}
\bibliography{bibfil}

\end{document}